\documentstyle[12pt]{article}

\newcommand{\const}{\mathop{\rm const}\limits}

\newcommand{\supp}{\mathop{\rm supp}\limits}

\begin{document}

\begin{center}

{\bf COMPOSED GRAND LEBESGUE SPACES} \\
\vspace{3mm}

 $ {\bf E.Ostrovsky^a, \ \ L.Sirota^b } $ \\

\vspace{4mm}

$ ^a $ Corresponding Author. Department of Mathematics and computer science, Bar-Ilan University, 84105, Ramat Gan, Israel.\\
\end{center}
E - mail: \ eugostrovsky@list.ru\\
\begin{center}
$ ^b $  Department of Mathematics and computer science. Bar-Ilan University,
84105, Ramat Gan, Israel.\\
\end{center}
E - mail: \ sirota@zahav.net.il\\

\vspace{3mm}

{\bf Abstract}.  In this article we introduce and investigate a new class of rearrangement invariant (r.i.) Banach function spaces, so-called Composed Grand Lebesgue Spaces (CGLS), in particular, Integral Grand Lebesgue Spaces (IGLS), which are some generalizations
of known Grand Lebesgue Spaces (GLS).\par
 We consider the fundamental functions of CGLS, calculate its Boyd's indices, obtain
 the norm boundedness some (regular and singular) operators in this spaces, investigate
 the conjugate and associate spaces, show that CGLS obeys the absolute continuous norm
 property etc.\par

\vspace{3mm}

\par

\vspace{3mm}

\vspace{3mm}

\section{Introduction. Definition of composed Grand Lebesgue Spaces.
Simple properties}

\vspace{3mm}

 {\bf 1.  Grand Lebesgue Spaces.} \par

 \vspace{2mm}

  Let $ (T, \cal{A},\mu) $ be some measurable space with sigma-finite
non-trivial measure $ \mu. $ For the measurable real valued
function $ f(t), \ t \in X, f: T \to R $ the symbol
$ |f|_p = |f|_p(X,\mu) $ will denote the usually $ L_p $ Lebesgue - M.Riesz
norm:
$$
|f|_p = ||f||L_p(X, \mu) = \left[ \int_X |f(t)|^p \ \mu(dt) \right]^{1/p}, \ p \ge 1.
\eqno(1.1)
$$
and correspondingly Lebesgue - M.Riesz spaces

$$
L_p = L_p(X,\mu) = \{ f: X \to R, |f|_p < \infty.  \}.
$$

\vspace{3mm}
We recall in this section  for reader convenience some definitions and facts from
the theory of (ordinary) Grand Lebesgue Spaces (GLS) spaces.\par
\vspace{3mm}

 Recently, see \cite{Fiorenza1},
 \cite{Fiorenza2}, \cite{Fiorenza3}, \cite{Iwaniec1}, \cite{Iwaniec2},
 \cite{Kozachenko1},\cite{Liflyand1}, \cite{Ostrovsky1}, \cite{Ostrovsky2},   etc.
 appears the so-called Grand Lebesgue Spaces $ GLS = G(\psi) =G\psi =
 G(\psi; a,b), \ a,b = \const, a \ge 1, a < b \le \infty, $ spaces consisting
 on all the measurable functions $ f: T \to R $ with finite norms

     $$
     ||f||G(\psi) \stackrel{def}{=} \sup_{p \in (a,b)} \left[ |f|_p /\psi(p) \right].
     \eqno(1.2)
     $$

     Here $ \psi(\cdot) $ is some continuous positive on the {\it open} interval
    $ (A,B) $ function such that

     $$
     \inf_{p \in (a,b)} \psi(p) > 0, \ \psi(p) = \infty, \ p \notin (a,b).
     $$

The set of all $ \psi $  functions with support $ \supp (\psi)= (a,b) $ will be
denoted by $ \Psi(a,b). $ \par
  This spaces are rearrangement invariant, see \cite{Bennet1}, chapter 2, and
  are used, for example, in the theory of probability  \cite{Kozachenko1},
  \cite{Ostrovsky1}, \cite{Ostrovsky2}; theory of Partial Differential Equations \cite{Fiorenza2}, 
  \cite{Iwaniec2};  functional analysis \cite{Fiorenza3},
  \cite{Iwaniec1},  \cite{Liflyand1}, \cite{Ostrovsky2}, \cite{Ostrovsky16}; theory of Fourier series \cite{Ostrovsky1}; theory of martingales \cite{Ostrovsky2}; mathematical statistics \cite{Sirota1}, 
  \cite{Sirota2}, \cite{Sirota4}, \cite{Sirota5}, \cite{Sirota6}, \cite{Sirota7}, \cite{Sirota8};  theory of approximation \cite{Ostrovsky7} etc.\par

\vspace{3mm}
{\bf  We introduce and investigate some extrapolation of GLS spaces,
a new class of rearrangement invariant (r.i.) Banach function
spaces, so-called Composed Grand Lebesgue Spaces (CGLS), in particular,
Integral Grand Lebesgue Spaces (IGLS). \par
  We consider the fundamental functions of CGLS, calculate its Boyd's indices, obtain
 the norm boundedness some (regular and singular) operators in this spaces, investigate
 the conjugate and associate spaces, show that CGLS obeys the absolute continuous norm
 property etc.} \par

\vspace{3mm}

{\sc We will denote by $  C_k = C_k(\cdot), k = 1,2,\ldots $ some positive finite essentially and by $ C, C_0 $ non-essentially "constructive" constants.}\par

\vspace{3mm}

{\bf 2. Composed Grand Lebesgue Spaces.}\par

\vspace{4mm}

  Let $ (X, \ ||\cdot||X) $ be a r.i. space, where $ X $ is linear subset on the space of all measurable function $ (a,b) \to R $ over our measurable space
 $ (T,M,\mu) $ with norm $ ||\cdot||X. $  \par

\vspace{3mm}

 {\bf Definition.} \par

\vspace{3mm}

 {\sc We will say that the space $ X $ with the norm $ ||\cdot||X $ is {\it Composed
Grand Lebesgue Space,} briefly: $ CGLS = CGLS(T,\mu; a,b; <>)=CGLS(T,\mu; <>)  $ space,
 if there exist a real constants $ a, b; 1 \le a < b \le \infty, $ and some {\it rearrangement invariant norm } $ < \cdot > $ defined on the space of a real functions defined on the interval $ (a,b), $ non necessary to be finite on all the functions such that}

  $$
  \forall f \in X \ \Rightarrow || f ||X = < \ h(\cdot) \ >, \ h(p) = |f|_p. \eqno(1.3)
  $$

\bigskip

{\it  Hereafter we will suppose that the norm $ ||\cdot||X $ is order continuous. }\par

 Another approach to the problem of syntheses of the spaces $ L_p $  see in 
\cite{Davis1}, \cite{Steigenwalt1}.\par

\bigskip

Recall (see \cite{Kantorovicz1}, chapter 4, section 3) that this means that for arbitrary sequence of a 
functions $ g_n(p), \ g_n(p) \in X $ the following implication holds:

$$
g_n(\cdot) \downarrow 0 \ \Longrightarrow \ <g_n> \downarrow 0. \eqno(1.4)
$$

 Note that the norm in the ordinary GLS is not order continuous. \par
\bigskip
    We will write for considered CGLS spaces $ (X, \ ||\cdot||X) $

   $$
       (a, b) \stackrel{def}{=} \supp(X),
   $$
   "moment support"; not necessary to be uniquely defined. But we will understand as
   the interval $ (a,b) $ its {\it minimal} value.\par

     It is obvious that arbitrary CGLS space is r.i. space.\par

   There are many r.i. spaces satisfied the condition (1.4): exponential Orlicz's spaces, 
some Marzinkievicz spaces, interpolation spaces (see [14], [42]).\par

{\it An important example.} \par
Let $ Q = \const \ge 1 $ and let $ \nu $ be some Borelian sigma-finite measure on the set 
$ (a,b). $  We  introduce the Integral Grand Lebesgue Space (IGLS)
$ G^{(Q)}_{\mu,\nu} = G^{(Q)}_{\mu,\nu}(a,b) $
as the set of all (measurable) functions $ f: T \in R $  with finite norm:

$$
||f||_Q^{(\mu,\nu)} \stackrel{def}{=} \left[\int_{[a,b]} |f|_p \ \nu(dp) \right]^{1/Q}.
\eqno(1.5)
$$

\vspace{3mm}
{\bf Remark 1.1.} If the measure $ \nu $ coincides with a Dirac  measure
$ \delta(p-p_0): $

$$
 \int_a^b h(p) \nu(dp) = f(p_0),
$$
then the Integral Grand Lebesgue Space (IGLS) $ G^{(Q)}_{\mu,\nu} $ is equal to the
Lebesgue-Riesz's space $ L_{p_0}. $\par

{\bf Remark 1.2.}
 Notice that if $ \nu \{ \infty \} = 0 $ or $ b < \infty, $ that this space satisfied the condition (1.4).\par
\vspace{3mm}
 {\it Subexample.}  In the article of S.F.Lukomsky \cite{Lukomsky1}
  are introduced the so-called $ G(p,\alpha) $ spaces
  consisted on all the measurable function $ f: T \to R $  with finite norm
   $$
   ||f||_{p,\alpha} = \left[ \int_1^{\infty} \left(\frac{|f|_x}{x^{\alpha}}
    \right)^p \ \nu(dx) \right]^{1/p},
   $$
   where $ \nu $ is some Borelian measure. Lukomsky  considered some applications of
   these spaces in the theory of Fourier series. \par
    Astashkin in \cite{Astashkin1}, \cite{Astashkin2}
    proved that the space $ G(p,\alpha) $ in the case
   $ T = [0,1] $ and $ \nu = m, \ m(\cdot) $ is Lebesgue measure, coincides with the
   Lorentz $ \Lambda_p( \log^{1-p \alpha}(2/s) ) $ space.  Therefore,
     both this spaces are CGLS spaces.\par
 Some  applications of the CGLS spaces in the approximation theory see in
 the article \cite{Sirota9}.\par

\vspace{3mm}
{\bf Theorem 1.1.} The Composed Grand Lebesgue Spaces (CGLS) are complete order continue
rearrangement invariant Banach functional spaces with Fatou property. \par
\vspace{3mm}
{\bf Proof.} It is suffices to prove only the completeness  of these spaces. In turn
it is suffices to prove the order continuity of the norm in this spaces, see, e.g.,
 \cite{Kantorovicz1}, chapter 4, section 3. \par
  Let $ g_n = g_n(t) $ be a monotonically decreasing sequence of a functions belonging to
the CGLS  such that as $ n \to \infty $ we have 

$$
g_n(t) \downarrow 0
$$
 $ \mu - $ almost everywhere.  As long as the classical Lebesgue-Riesz's space $ L_p $
 obeys the order continuity of the norm, as $ n \to \infty $

$$
|g_n(\cdot)|_p \downarrow 0
$$
for all the values $ p; \ p \in (a,b). $.\par
 Since the norm $ <\cdot> $ has the order continuity of the norm, we conclude

 $$
 ||g_n||X = < |g_n(\cdot)|_p > \downarrow 0.
 $$
This completes the proof of theorem 1.1.\par
 For instance, the last passing to the limit as $ n \to \infty $ for the Integral Grand
Lebesgue Spaces may be grounded by means of Lebesgue dominated convergence theorem.\par

\vspace{3mm}

\section{ Norm absolutely continuity (ACN)}

\vspace{3mm}

 Recall that the function $ f $ from the rearrangement invariant space $  X $ with the norm $ ||\cdot||X $ has absolutely continuous norm: $ f \in ACN $ in the terminology, e.g. of the book \cite{Bennet1}, chapters 2,3,  if for arbitrary
 decreasing sequence of measurable sets $ E_n $ such that $ \mu -$ almost everywhere

 $$
 \cap_n E_n = \emptyset
 $$
there holds

$$
\lim_{n \to \infty} ||f \cdot I(E_n)||X = 0.
$$
Hereafter $ I(A) = I(A,t) $ denotes the indicator function for the set $ A:$

$$
I(A,t) = 1, \ t \in A; \ I(A,t) = 0, \ t \notin A.
$$
 The whole space $ X $ has ACN property, iff arbitrary function $ f,  f \in X $ has
 ACN property. \par

\vspace{2mm}
{\bf Theorem 2.1} The CGLS space $ X = CGLS(T,\mu; a,b; <>) $ has ACN property.\par
{\bf Proof} is at the same as in theorem 1.1. Indeed, let

 $$
  E_{n+1} \subset E_n, \  \cap_n E_n = \emptyset,
 $$
and let $ f \in X; $ then $ |f|\cdot I(E_n) \in X $  and
$ |f|\cdot I(E_n) \downarrow 0 \ \mu - $ a.e.  Therefore \\
$ |||f|\cdot I(E_n)||X \downarrow 0, $ Q.E.D.\par

\vspace{2mm}

{\bf Corollary 2.1.} The CGLS space $ X = CGLS(T,\mu; a,b; <>) $ has {\it Lebesgue,} or
{\it dominated convergence } property:

$$
\forall \{ f_n \}: |f_n| \le f, \ f,f_n \in X, \
f_n \stackrel{\mu - a.e.}{\to} g \ \Rightarrow ||f_n - g||X \to 0.
$$

\vspace{3mm}

\section{Separability and reflexivity.}\par

\vspace{3mm}

{\bf Theorem 3.1.} Assume that the measure $  \mu $ on the set $ A $
is separable; this imply by definition that the space $ \cal{A} $ is separable relatively a distance

$$
d(B_1,B_2) \stackrel{def}{=} \mu(B_1 \setminus B_2) + \mu(B_2 \setminus B_1)=
\mu(B_1 \Delta B_2), \ B_1,B_2 \in \cal{A}.
$$
 Then the space $ CGLS(T,\mu; <>)$ is separable.\par
 {\bf Proof } follows immediately from the proved order continuity of this space and
 from the Theorem 3 of monograph Kantorovicz L.V., Akilov G.P.
 \cite{Kantorovicz1}, chapter 4,  section 3. \par

\vspace{2mm}

{\bf Theorem 3.2.} Consider the space (IGLS) $ G^{(Q)}_{\mu,\nu} $  on the set
$ p \in [a,b], $ where $ \nu{\{1 \}} = \nu{\{\infty \}} = 0 $ or $ 1 < a < b < \infty. $\par
 If the measure $ \nu $ is purely atomic:

 $$
 \int_a^b h(p) \ \nu(dp) = \sum_{k=1}^n c(k) f(p_k), \ c(k) > 0, \ p_k \in (a,b),
 $$
then the space  $ G^{(Q)}_{\mu,\nu} $ is reflexive.\par
{\bf Proof} is obvious:

$$
||f||G^{(Q)}_{\mu,\nu} = \left[\sum_{k=1}^n c(k)|f|_{p_k}^Q \right]^{1/Q}
$$
and follows from the reflexivity of all the  spaces $ L_{p_k}. $ Note also that 
in considered case the space $ G^{(Q)}_{\mu,\nu} $ coincides up to norm equivalence 
with  the direct sum of $ L_{p_k} $ spaces. \par

\vspace{3mm}

{\bf Remark 3.1.} The author does not know  an essential generalization of this result.\par

\vspace{3mm}
\section{ Boyd's indices  for IGLS spaces}
\vspace{3mm}

 Let in this section $ T = R_+ = (0,\infty) $ with Lebesgue measure. Introduce the so-called 
 dilation operator (more exactly, the family of operators) $ \sigma_s $ by the formula

 $$
 \sigma_s [f](t) = f(t/s).\eqno(4.1)
 $$

 For arbitrary r.i. space over $ T= R_+, $ for example, for the space
  $ Y=CGLS(T,\mu; a,b; <> $ the Boyd's indices $ \alpha(Y), \beta(Y)  $ may be defined
  as follows:

  $$
  \beta(Y) \stackrel{def}{=} \lim_{s \to \infty} \frac{\log ||\sigma_s||Y}{|\log s|},
  $$

  $$
  \alpha(Y) \stackrel{def}{=} \lim_{s \to 0+} \frac{\log ||\sigma_s||Y}{|\log s|}.
  $$
These limits there exists and play a very important role in the theory of Fourier series
and in the theory of singular integral operators in r.i. spaces, see \cite{Bennet1},
chapters 5,6.\par

{\bf Theorem 4.1.} Let $ \nu $ be Lebesgue measure:
$$
\nu(A) = \int_A  dp,
$$

 We assert for the IGLS space $ G^{(Q)}_{\mu,\nu; a,b}: $

$$
\alpha(G^{(Q)}_{\mu,\nu; a,b}) = \frac{1}{b},\eqno(4.2)
$$

$$
\beta(G^{(Q)}_{\mu,\nu; a,b}) = \frac{1}{a}, \eqno(4.3)
$$

{\bf Proof.} We can suppose without loss of generality that $ b < \infty $
and that $ \nu \{[a,b]  \} = 1. $\par
It is sufficient also to investigate only upper Boyd's index
$ \beta(G^{(Q)}_{\mu,\nu; a,b}); $ the case of lower index
$ \alpha(G^{(Q)}_{\mu,\nu; a,b}) $ is investigated analogously. \par
{\bf A. Upper bound.} Let $ f = f(t) $ be arbitrary positive function from the space
$ G^{(Q)}_{\mu,\nu; a,b}$. We have consequently:

$$
|\sigma_s [f]|_p^p = \int_0^{\infty} |f(t/s)|^p dt =  \int_0^{\infty} s
\int_0^{\infty} |f(t)|^p dt = s |f|_p^p;
$$

$$
|\sigma_s [f]|_p =  \le s^{1/a} |f|_p, \ s > 2;
$$

$$
||\sigma_s|| \le s^{1/a}, \ \beta(G^{(Q)}_{\mu,\nu; a,b})\le 1/a.
$$

{\bf B. Lover bound.} We set for simplicity $ Q=1. $ We have for at the same continuous
positive function $ f$  from the space $ \ G^{(Q)}_{\mu,\nu; a,b} \ $ and each "small"
positive number $ \epsilon:  $

$$
I:= \int_a^b s^{1/p} |f|_p dp \ge s^{1/(a+\epsilon)} \int_a^{a+\epsilon}|f|_p dp, \
s > 2;
$$
Therefore

$$
\frac{||\sigma_s||}{\log s} \ge \frac{1}{a+\epsilon} - \frac{C(\epsilon)}{|\log s|};
$$

$$
\overline{\lim}_{s \to \infty} \frac{||\sigma_s||}{\log s} \ge \frac{1}{a}.
$$

{\bf Remark 4.1.} At the same result with at the same proof is true for arbitrary
$ Y=CGLS(T,\mu; a,b; <> $ space.\par

 But for the $ G^{(Q)}_{\mu,\nu; a,b} $ space with $ T=R_+ $ for the Boyd's index
 our result may be refined. Namely, let the function $ f, f:T \in R $ be such that
 the function $ p \to |f|_p $ is positive, bilateral bounded  and is  continuous in
 some neighborhood of the point $ p=a: p \in (a,a+\epsilon).$ Then we conclude
 as $ s \to \infty, s > 3 $ using the saddle-point method:

 $$
 ||\sigma_s||^Q \ge C \int_a^b s^{Q/p} |f|_p^Q dp \ge C_2 \int_a^b s^{Q/p} dp \ge
C_3 C_4^Q \frac{s^{Q/a}}{\log s};
$$

$$
||\sigma_s|| \ge C_5 \frac{s^{1/a}}{|\log s|^{1/Q}}.
$$
 Therefore, we conclude for sufficiently greatest values $  s: $

$$
\frac{1}{a} - \frac{\log |\log s|}{Q |\log s|} - \frac{C_6}{|\log s| }
\le \frac{\log ||\sigma_s||}{|\log s|} \le \frac{1}{a} \eqno(4.4)
$$
and we find analogously in the case when $ b < \infty $ for the smallest values
$s, \ s \in (0, e^{-2e}):  $

$$
\frac{1}{b} \le \frac{\log ||\sigma_s||}{|\log s|} \le \frac{1}{b} +
\frac{\log |\log s|}{Q |\log s|} + \frac{C_7}{|\log s| }.\eqno(4.5)
$$

\vspace{2mm}
\section{Fundamental function}
\vspace{2mm}

 Recall that the fundamental function $ \phi = \phi_Y(\delta), \ \delta \in (0,\infty) $
 for arbitrary r.i.  space $ Y $  with norm $ ||\cdot||Y $ is defined by formula

$$
\phi_Y(\delta) = || I(A)||Y. \eqno(5.1)
$$

{\it In this section we investigate some properties of the fundamental function for the
 Integral Grand Lebesgue Spaces } (IGLS) $ Y= G^{(Q)}_{\mu,\nu;a,b}. $ \par

 Note that in the considered case

 $$
 \phi_Y^Q(\delta) = \phi_{G^{(Q)}_{\mu,\nu;[a,b]}}^Q(\delta) =
 \int_{[a,b]} \delta^{Q/p} \nu(dp), \eqno(5.2)
 $$
or after substitution $ \delta = \exp(\lambda/Q), \ \lambda \in (-\infty,\infty): $

$$
\zeta_Y(\lambda) =
\zeta(\lambda) = \zeta_Y^{(Q)}(\lambda) \stackrel{def}{=} \phi_Y^Q
\left(e^{\lambda/Q} \right)=
\int_{[a,b]} e^{\lambda/p} \nu(dp). \eqno(5.3)
$$

{\bf Theorem 5.1.} Let $ \nu([a,b]) < \infty. $  \par
{\bf A.}  The function $ \zeta_Y(\lambda) $ is absolutely monotonic on the set
$ \lambda \in (-\infty,\infty):  $

$$
\zeta^{(k)}(\lambda) > 0, \ k = 0,1,2,\ldots \eqno(5.4)
$$
and such that

$$
C_1 e^{\lambda/b} \le \zeta(\lambda) \le C_2 e^{\lambda/a}. \eqno(5.5)
$$

  We denote further the class of a functions satisfying the conditions (5.4) and (5.5)
  as $ AM = AM(a,b). $ For instance, any function from the class  $ AM(a,b) $ is infinitely
  differentiable on the whole axis $ R^1. $\par
\vspace{2mm}
{\bf B.} Conversely, if the function $ \lambda \to \zeta(\lambda), $ defined on the set
$ \lambda \in (-\infty,\infty) $ satisfies the conditions (5.4) and (5.5), then
there exists a finite Borelian measure $  $ on the closed interval $ [a,b] $ for which

$$
\zeta(\lambda) =
\int_{[a,b]} e^{\lambda/p} \nu(dp). \eqno(5.6)
$$
{\bf Proof.} The proof of the assertion {\bf A } is very simple. Indeed:

$$
 \zeta^{(k)}(\lambda) = \int_{[a,b]} p^{-k} \ e^{\lambda/p} \nu(dp) > 0;
$$

$$
\zeta(\lambda) \le e^{\lambda/a} \ \nu([a,b]) = C_2 e^{\lambda/a},
$$

$$
\zeta(\lambda) \ge e^{\lambda/b} \ \nu([a,b]) = C_1 e^{\lambda/b}.
$$

{\bf Proof} of the assertion {\bf B. }  Let the inequality (5.4) holds. From the theorem
of S.N. Bernstein follows that there exists a finite measure $ \tilde{\nu} $ on the real
axis $ (-\infty,\infty) $ for which

$$
\zeta(\lambda) = \int_{[0,\infty)}e^{\lambda x} \tilde{\nu}(dx).
$$

We conclude by virtue of the inequalities (5.5) that the support of the measure
$  \tilde{\nu} $ contained in the set $ x \in [1/b, 1/a]: $

$$
\zeta(\lambda) = \int_{[1/b,1/a]} e^{\lambda x} \tilde{\nu}(dx).\eqno(5.7)
$$
 The assertion {\bf B} of theorem (5.1) follows  from (5.7) after substitution $ x=1/p. $
 \par
{\bf Remark 5.1.} Note in additional that the function $ \zeta(\lambda) $ satisfies the
following  inequalities:

$$
b^{-k} \zeta^{(k)}(\lambda)\le \zeta(\lambda) \le a^{-k} \ \zeta^{(k)}(\lambda), \
\lambda \in R, k=0,1,2,\ldots.
$$
\vspace{2mm}
 We continue the investigation of the properties of the fundamental function for the
spaces $ Y= G^{(Q)}_{\mu,\nu;a,b}, \ a < b < \infty. $  \par

\vspace{3mm}
{\bf Theorem 5.2.}
Suppose the measure $ \nu(\cdot) $ is absolutely continuous relative Lebesgue measure:

$$
\nu(B) = \int_B h(p) dp,
$$
where the function $ h(p)$ is non-negative, integrable and continuous at the points
$ p=a $ and $ p = b $ such that $ h(a) > 0, h(b) > 0. $\par
 Then as $ \delta \to \infty $

$$
\phi_Y(\delta) \sim \left[ \frac{h(a)}{Q \log \delta} \right]^{1/Q} \cdot
\delta^{1/a},\eqno(5.8)
$$

and as $ \delta \to 0+ $

$$
\phi_Y(\delta) \sim \left[ \frac{h(b)}{Q \log \delta} \right]^{1/Q} \cdot
\delta^{1/b}. \eqno(5.9)
$$
{\bf Proof } follows immediately from the theory of saddle-point method, see
\cite{Fedoruk1}, chapter 2, section 4.\par

\vspace{3mm}

{\bf Example 5.1.} Let

$$
\nu(B) = \int_B p^{-2} dp,
$$
then
$$
\phi_Y^Q(\delta) = \int_a^b \delta^{Q/p} \frac{dp}{p^2} =
\frac{\delta^{Q/a}-\delta^{Q/b}}{Q \log \delta},
$$
therefore

$$
\phi_Y(\delta) = \left[\frac{\delta^{Q/a}-\delta^{Q/b}}{Q \log \delta}\right]^{1/Q}.
$$

\vspace{3mm}

\section{ Conjugate and associate spaces}

\vspace{3mm}

  Since the CGLS space $ Y $ has the ACN property, the conjugate space $ Y^*$ coincides 
with associate space $ Y'.$ Therefore. every continuous (bounded) linear functional
 $ l, l: Y \to R $ has a view

 $$
 l(f) = l_g(f) = \int_T f(t) g(t) \mu(dt),\eqno(6.1)
 $$
 where $ g:T \to R $ is some measurable function for which

 $$
 ||g||Y^* = ||g||Y' = \sup_{f: ||f||Y = 1} \left|\int_T f(t) g(t) \mu(dt) \right| < \infty.\eqno(6.2)
 $$

  But the expression (6.2) is very hard to calculate. We give now a simple upper
  estimation for the norm $ ||g||Y^* .$ \par
 Let us consider the Integral Grand Lebesgue Space (IGLS) $ G^{(Q)}_{\mu,\nu}. $ We
 denote as usually $ Q'= Q/(Q-1), \ p'= p/(p-1), \  Q,p > 1 $ and define

 $$
 ||g||_{Q'} = ||g||_{Q',\nu} = \left[ \int_a^b |g|_{p'}^{Q'} \ \nu(dp)  \right]^{1/Q'}.
 $$

 {\bf Theorem 6.1.} We assert for the space $ Y =G^{(Q)}_{\mu,\nu}: $

 $$
 ||g||Y^* \le ||g||_{Q'}. \eqno (6.3)
 $$

 {\bf Proof.}  As long as the measure $ \nu $ is sigma-finite, we can and will suppose
 $ \nu((a,b))=1. $ We obtain using H\"older's inequality for the representation (6.1):

 $$
 |l_g(f)| = \left| \int_T f(t) g(t) \mu(dt) \right| \le |f|_p |g|_{p'};
 $$
therefore
 $$
|l_g(f)| \le \inf_{p \in (a,b)} \left[ |f|_p |g|_{p'} \right] \le
\int_a^b |f|_p |g|_{p'} \nu(dp) \le
 $$

 $$
 \left[\int_a^b |f|_p^Q \nu(dp) \right]^{1/Q} \cdot
     \left[\int_a^b |g|_{p'}^{Q'} \nu(dp)  \right]^{1/Q'} = ||f||Y \cdot ||g||_{Q'},
 $$
 Q.E.D.\par
 Note that the estimation (6.3) of theorem 6.1 is exact when  the measure $ \nu $ is
 Dirac measure  and is exact up to multiplicative constant for the pure discrete
 measure $  \nu $ with finite support.\par
  But in general case the expression for $ ||g||_{Q'}$   does not represent the general
 form for linear functional over the space $ G^{(Q)}_{\mu,\nu} = G^{(Q)}_{\mu,\nu}(a,b).$
 Let us consider correspondent example. \par
\vspace{3mm}
{\bf Theorem 6.2.}
  Suppose the measure $ \nu $ is absolutely continuous over Lebesgue measure with positive continuous differentiable Radon-Nikodim derivative $ h(p): $

  $$
  \nu(B) = \int_B h(p) dp.
  $$
where the function $ h(p)$ is non-negative, integrable and continuous at the points
$ p=a $ {\it or} at the point $ p = b $ and such that $ h(a) > 0, h(b) > 0. $ Then the
Integral Grand Lebesgue Space (IGLS) $ G^{(Q')}_{\mu,\nu}( b/(b-1), a/(a-1)), a > 1,
b < \infty $  does not coincide with the dual (associate) space to the space
$ G^{(Q)}_{\mu,\nu}(a, b).  $\par
{\bf Proof.} Assume conversely. i.e. that

$$
G^{(Q')}_{\mu,\nu} (b/(b-1), a/(a-1)) = G^{(Q)}_{\mu,\nu} (a, b).
$$
 It is known (see \cite{Bennet1}), chapter 1,2 that

 $$
 \phi_Y(\delta)\cdot \phi_{Y'}(\delta) = \delta, \ \delta \in (0,\infty).
 $$

 The last equality contradict the proposition of theorem 5.2 when $ \delta \to \infty $ or
as $ \delta \to 0+. $\par

\vspace{3mm}
\section{Relations with another r.i. spaces.}
\vspace{3mm}

{\it  We intend to prove in this section that the classical r.i. spaces: Lorentz,
Orlicz, Marzinkievicz, Grand Lebesgue Spaces  does not coincide or equivalent to the
considered in this article Composed Grand Lebesgue Spaces.} \par

 Note first of all that the {\it exponential} Orlicz's and Grand Lebesgue Spaces are under simple conditions not separable \cite{Liflyand1}, \cite{Ostrovsky19}, in contradiction to the properties of CGLS spaces.\par
  We introduce {\it in this section} the following equivalence relation $ \asymp $
  between a two positive functions $ g_1(\lambda), g_2(\lambda) $ defined on the whole
  real axis:

  $$
  g_1(\cdot) \asymp g_2(\cdot) \Leftrightarrow \exists C_1,C_2 = \const >0,  C_1 \le C_2, \  C_1 \le \frac{g_1(\lambda)}{g_2(\lambda)} \le  C_2. \eqno(7.0)
  $$

  Now, let $ X, \ ||\cdot||X $ be r.i. space over our measurable triple
 $ (T, \cal{A},\mu); $ denote by $ \phi_X(\delta) $ its fundamental function and define
  for some $ Q = \const \ge 1 $ and for arbitrary  $ \lambda \in R $ the following
 function:

$$
\tau_X(\lambda) = \tau_X^{(Q)}(\lambda) := \phi_X^Q\left(e^{\lambda/Q} \right).
\eqno(7.1)
$$
 It  follows from the properties of the fundamental functions of CGLS spaces
 (theorem 5.1) the following result.\par
{\bf Theorem 7.1.} \par
{\bf A.} If

$$
\tau_X^{(Q)}(\cdot) \notin \cup_{(a,b): \ 1 \le a < b < \infty} \{ AM(a,b) \}, \eqno(7.2)
$$
then the r.i. space $ X $ does not coincide with arbitrary $ G^{(Q)}_{\mu,\nu}(a,b) $
space.\par
{\bf 2.} If for fixed r.i. space $ X $ and for all the values $ (a,b), \ Q $
its function $\tau_X^{(Q)}(\cdot) $ is not equivalent in the relation $ \asymp $
to arbitrary function from the class $ AM(a,b), $ then the space $ X $ is not equivalent
(in the sense of Banach space equivalency) to arbitrary $ G^{(Q)}_{\mu,\nu}(a,b) $ space.\par

 For instance, the classical r.i. Lorentz and Marzinkievicz spaces, see \cite{Krein1}, chapter3.
may have the non-smooth fundamental function; therefore this spaces does not coincide with 
arbitrary $ G^{(Q)}_{\mu,\nu}(a,b) $ space.\par

\vspace{3mm}
\section{Convergence and compactness}
\vspace{3mm}

 As we know, the space $ G^{(Q)}_{\mu,\nu}(a,b) $ obeys the Absolutely Continuous Norm
 property. We conclude consequently: \par

 \vspace{3mm}

 {\bf Theorem 8.1.} Let $  F = \{ f_{\alpha} \} \subset G^{(Q)}_{\mu,\nu}(a,b),  \alpha \in A,  $
 where $ A $ is arbitrary set of indices be any subset of the space $ G^{(Q)}_{\mu,\nu}(a,b).$ \par
  The set $  F  $ is compact set in this space iff it is bounded, closed and
obeys the Uniform Absolutely Continuous Norm property. \par

\vspace{3mm}
{\bf Proof} follows immediately from the theory of a rearrangement invariant spaces, see
\cite{Bennet1}, chapter 2; \cite{Kantorovicz1}, chapter 4, section 3. \par
 A slight consequence: \par
\vspace{3mm}
{\bf Theorem 8.2.} The sequence $ \{f_n \}, n-1,2,\ldots $ of a functions from the space
$ G^{(Q)}_{\mu,\nu}(a,b)$ converges as $ n \to \infty $ iff it converges in $ \mu \ - $
measure and has the Uniform Absolutely Continuous Norm property. \par
\vspace{3mm}
 Obviously, the convergence in measure may be replaced by convergence in some
 $  L_p, \ p \in (a,b).  $\par

  Let $ (T, \cal{A},\mu) $ be again some measurable space with sigma-finite
non - trivial measure $ \mu. $ We assume in addition that $ T = \{t\} $ is homogenous
compact metric space with additive operation $ t \pm s. $   We define the difference 
operator $ \Theta_h[f], \ h,t \in T $ as ordinary:

$$
\Theta_h[f](t) = f(t+h) - f(t).
$$

The application of Shilov's theorem give us the following result: \par

\vspace{3mm}

 {\bf Theorem 8.3.} Let $  F = \{ f_{\alpha} \} \subset G^{(Q)}_{\mu,\nu}(a,b),  \alpha \in A,  $
 where $ A $ is arbitrary set of indices be any subset of the space $ G^{(Q)}_{\mu,\nu}(a,b).$ \par
  The set $  F  $ is compact set in this space iff it is bounded, closed and

$$
\lim_{h \to 0} \sup_{f \in F} <\Theta_h \left[f_{\alpha} \right](\cdot) > = 0.
$$

\vspace{3mm}
\section{Some applications: boundedness of Hilbert's transform  and other operators in CGLS}
\vspace{3mm}
  Let $ 1 \le a < b \le \infty $ and let $ U $ be an operator, not necessary to be linear or sublinear,
such that there exists non-zero finite "constant" $ K = K(p), \ a < p < b $ for which

$$
|U[f] |L_p(X, \mu)  \le K(p) \ |f|L_p(X, \mu), \ f \in L_p(X, \mu). \eqno(9.1)
$$
 The condition (9.1) is satisfied, e.g. for weight Fourier transform \cite{Beckner1}, \cite{Leindler1},
 for Hilbert's transform: $ a = 1, b = \infty, \ K(p) = Cp^2/(p-1) $ \cite{Pichorides1},
 Hardy-Littlewood maximal operator:
  $ a = 1, b = \infty, \ K(p) = Cp/(p-1) $ \cite{Hardy1},  mean value operator

 $$
 U[a](n) = n^{-1} \sum_{k=1}^n a(k), \ K(p) = Cp/(p-1),
 $$
 \cite{Hardy1}, and many others singular operators.\par
We will understood furthermore as the value $ K(p) $ its minimal value, namely:

$$
K(p) = \sup_{f \in L_p, f \ne 0} \left[ \frac{|U[f] |L_p(X, \mu)}{|f|L_p(X, \mu)} \right]. \eqno(9.2)
$$
 This constants is calculated in many works, see e.g. in \cite{Hardy1}, \cite{Leindler1},
 \cite{Pichorides1}, \cite{Talenti1}.\par
We introduce a new CGLS space as a set of all (measurable) functions with finite norm

$$
|||g||| \stackrel{def}{=} < |g|_p/K(p) >. \eqno(9.3)
$$

{\bf Theorem 9.1.}

$$
||| U[f]||| \le 1 \cdot <f>, \eqno(9.4)
$$
where the constant "1" in (9.4) is the best possible.\par
 {\bf Proof.} The lower bound is attained, for instance, when the measure
 $ \nu $ is Dirac's $ \delta(p-p_0) $ measure. The upper bound may be proved 
 very simple:

 $$
 ||| U[f] ||| = \ < |U[f]|_p/K(p)> \ \le \ < K(p) |f|_p/K(p) > \ =
 $$
 $$
  \ < |f|_p > \ = \ < f >. \eqno(9.5)
 $$

\end{document}